\documentclass{amsart}

\usepackage{amsfonts,amssymb,verbatim,amsmath,amsthm,latexsym,textcomp,amscd}
\usepackage{latexsym,amsfonts,amssymb,epsfig,verbatim}
\usepackage{amsmath,amsthm,amssymb,latexsym,graphics,textcomp}
\usepackage{graphicx}
\usepackage{color}
\usepackage{url}

\input xy
\xyoption{all}

\theoremstyle{plain}
\newtheorem{theorem}{Theorem}[section]

\newtheorem{question}[theorem]{Question}

\newtheorem{obs}[theorem]{Observation}

\newtheorem{prop}[theorem]{Proposition}
\newtheorem{lemma}[theorem]{Lemma}
\newtheorem{cor}[theorem]{Corollary}
\newtheorem{defn}[theorem]{Definition}

\newtheorem{claim}[theorem]{Claim}
\newtheorem{rmk}[theorem]{Remark}

\newcommand{\bdy}{\partial}

\newcommand{\E}{{\mathbb E}}

\newcommand{\natls}{{\mathbb N}}
\newcommand{\ratls}{{\mathbb Q}}

\newcommand\Hyp{{\mathbf H}}
\newcommand\D{{\mathbf D}}

\newcommand\CC{{\mathcal C}}

\newcommand\EE{{\mathcal E}}
\newcommand\FF{{\mathcal F}}
\newcommand\GG{{\mathcal G}}

\newcommand\JJ{{\mathcal J}}

\newcommand\LL{{\mathcal L}}
\newcommand\MM{{\mathcal M}}

\newcommand\PP{{\mathcal P}}

\newcommand\VV{{\mathcal V}}

\newcommand\PMF{{\PP\kern-2pt\MM\FF}}

\newcommand\PML{{\PP\kern-2pt\MM\LL}}

\newcommand{\fsubd}{\mathrel{{\scriptstyle\searrow}\kern-1ex^d\kern0.5ex}}
\newcommand{\bsubd}{\mathrel{{\scriptstyle\swarrow}\kern-1.6ex^d\kern0.8ex}}
\newcommand{\fsubeq}{\mathrel{\raise-.7ex\Hbox{$\overset{\searrow}{=}$}}}
\newcommand{\bsubeq}{\mathrel{\raise-.7ex\Hbox{$\overset{\swarrow}{=}$}}}

\newcommand{\co}{{\colon}}

\begin{document}

\title{Pattern Rigidity in Hyperbolic Spaces:\\ Duality
and PD Subgroups}
\author{Kingshook Biswas}
\author{Mahan Mj}

\thanks{Research partly supported by a DST Project grant \\ AMS subject classification =   20F67(Primary), 22E40   57M50}

\begin{abstract} 
For $i= 1,2$, let $G_i$ be   cocompact groups
of isometries of hyperbolic space  
 $\Hyp^n$ of real dimension $n$, $n \geq 3$. 
Let $H_i \subset G_i$ be 
 infinite index quasiconvex subgroups satisfying one of the following
conditions: \\
1) limit set of $H_i$ is a codimension
one topological sphere.\\
2) limit set of $H_i$ is an even dimensional topological sphere.\\ 
3) $H_i$ is a codimension one {\it duality group}. This generalizes 
(1). In particular, if  $n = 3$, $H_i$ could be {\it any} freely 
indecomposable subgroup of $G_i$. \\
4) $H_i$ is an odd-dimensional  Poincar\'e Duality group $PD(2k+1)$.
This generalizes (2).\\
  We prove pattern rigidity for such
pairs extending work of Schwartz who proved pattern rigidity when $H_i$
is cyclic. All this generalizes to quasiconvex subgroups
of uniform lattices in rank one symmetric spaces satisfying one of
the conditions (1)-(4), as well as certain
special subgroups with disconnected limit sets.
In particular, pattern rigidity holds for {\bf all quasiconvex 
subgroups of hyperbolic 3-manifolds that are not virtually
free}.
Combining this with the main result of Mosher-Sageev-Whyte
\cite{msw2}, we get quasi-isometric rigidity results for graphs
of groups where the vertex groups are uniform 
lattices in rank one symmetric spaces and edge groups are of any
of the above types. \\

\medskip

\begin{center}
{\em This paper is dedicated to the memory of Kalyan Mukherjea.}
\end{center}

\end{abstract}

\maketitle

\tableofcontents

\section{Introduction}
\subsection{Statement of Results}
In \cite{gromov-ai} Gromov proposed the project of classifying
 finitely generated groups up to quasi-isometry. A class of
groups where
any two members are quasi-isometric if and only if they are commensurable
is said to be quasi-isometrically rigid.
However, in certain
classes
of groups, for instance uniform lattices $G$ in some fixed hyperbolic 
space $\Hyp$,  all members of the class are
quasi-isometric to $\Hyp$ and hence
to each other. In this context (or in a context where quasi-isometric
rigidity is not known) it makes sense to ask a relative version
of Gromov's question. Here, (almost as a rule) additional restriction
is imposed on the quasi-isometries by requiring that they preserve some 
additional structure given by a `symmetric pattern' of subsets. A 
`symmetric pattern' of subsets
 roughly means a $G$-equivariant  collection $\JJ$ of subsets
in $\Hyp$, each of which in turn is invariant under a conjugate
of a fixed subgroup $H$ of $G$, such that the quotient of an element
of $\JJ$ by its stabilizer is compact. Then the relative version
of Gromov's question for classes of pairs $(G,H)$ becomes: \\

\begin{question} Given a quasi-isometry $q$ of two such pairs
$(G_i, H_i)$ ($i = 1, 2$) pairing a $(G_1, H_1)$-symmetric pattern
$\JJ_1$ with a $(G_2, H_2)$-symmetric pattern
$\JJ_2$, does there exist an isometry $I$ which performs the
same pairing? Further, does $q$ lie within a bounded distance
of $I$ ?
\label{motvn}
\end{question}

Formulated in these terms, the phenomenon addressed
by Question \ref{motvn} is called {\it pattern rigidity}
(See \cite{msw2} where this terminology was first used. See Section
\ref{outline} for more on the genesis of the problem and the
techniques used, particularly  work of Mostow
\cite{mostow-pihes} and Sullivan \cite{sullivan-rigidity}.)

One of
the first papers to come out in the subject of quasi-isometric rigidity
was by Schwartz \cite{schwartz_pihes}, and even here, the problem can be formulated (in part) as a pattern rigidity question for symmetric
patterns 
of horoballs in $\Hyp$. The next major piece of work on pattern rigidity
was for subgroups $H = \mathbb{Z}$ by Schwartz \cite{schwarz-inv}
again. In a certain sense, \cite{schwartz_pihes} deals with 
symmetric
patterns of convex sets whose limit sets are single points, and 
\cite{schwarz-inv} deals with 
symmetric
patterns of convex sets (geodesics)
whose limit sets consist of two points.  In this paper we
initiate the study of pattern rigidity for
symmetric
patterns of convex sets whose limit sets are infinite. 

For $i= 1,2$, let $G_i$ be   cocompact groups
of isometries in a rank one symmetric space $\Hyp^n$ of real dimension $n$, $n \geq 3$. 
Let $H_i \subset G_i$ be an
 infinite index quasiconvex subgroup satisfying {\bf any} of the following
conditions: 
\begin{enumerate}
\item limit set of $H_i$ is a codimension
one topological sphere.
\item  limit set of $H_i$ is an even dimensional topological sphere.
\item  $H_i$ is a codimension one {\it duality group}. This generalizes 
(1).  
\item  $H_i$ is an odd-dimensional  Poincar\'e Duality group $PD(2k+1)$.
This generalizes (2).
\end{enumerate}

In this paper, 
we prove pattern rigidity for such
pairs (See Theorem \ref{omni} below for a precise statement).

Examples of (1) above include quasi-Fuchsian surface subgroups of closed hyperbolic 3-manifold groups
corresponding to immersed surfaces. A special case of this would correspond to {\it embedded totally geodesic}
surfaces (cf. \cite{frigerio}). Examples of (3) include all freely 
indecomposable subgroup of closed hyperbolic 3-manifold groups. (Note that these need not correspond 
to embedded submanifolds with boundary.) See also \cite{belegradek-hyperplane}, \cite{bi} and \cite{mahan-hs}
for related work.

\begin{defn}
A {\bf symmetric pattern}
 of closed convex (or quasiconvex) sets in a rank one symmetric
space $\Hyp$ is a $G$-invariant countable collection 
$\JJ$ of closed convex (or quasiconvex) sets such that 

\begin{enumerate}
\item $G$ acts cocompactly on $\Hyp$. 
\item The stabilizer $H$ of $J \in \JJ$ acts cocompactly on $J$. 
\item  $\JJ$ is the orbit of some (any) $J \in \JJ$ under $G$. 
\end{enumerate}
\end{defn}

This definition is slightly more restrictive than Schwartz' notion
of a symmetric pattern
 of geodesics, in the sense that he takes
$\JJ$ to be a finite union of orbits of geodesics, whereas Condition (3)
above forces $\JJ$ to consist of one orbit.
All our results go through with the more general definition, where
$\JJ$ is a finite union of orbits of closed convex (or quasiconvex) sets,
but we restrict ourselves to one orbit for expository ease.

Suppose that
$(X_1, d_1), (X_2,d_2)$ are metric spaces. Let $\JJ_1, \JJ_2$ be collections of closed subsets of $X_1, X_2$ respectively. Then $d_i$ induces a non-negative function (which, by abuse of notation, we continue to refer to as $d_i$) on $\JJ_i \times \JJ_i$ 
for $i = 1, 2$. This is just the ordinary (not Hausdorff) distance between closed subsets of a metric space. 

In particular, 
consider two hyperbolic groups $G_1, G_2$ with quasiconvex subgroups $H_1, H_2$, Cayley graphs $\Gamma_1, \Gamma_2$. Let $\LL_j$  for $j = 1, 2$ denote the collection of translates of limit sets of $H_1, H_2$ in $\partial G_1, \partial G_2$ respectively. Individual members of  the collection $\LL_j$ will be denoted as $L^j_i$.
Let $\JJ_j$ denote the collection 
$\{ J_i^j = J(L_i^j): L_i^j \in \LL_j \}$ of joins of limit sets.
Recall that the join of a limit set
$\Lambda_i$ is the union of bi-infinite geodesics in $\Gamma_i$ with end-points in $\Lambda_i$. This is a uniformly quasiconvex set and lies at a bounded Hausdorff distance from the Cayley graph of the subgroup $H_i$ 
Following Schwartz \cite{schwarz-inv}, we define:

\begin{defn} A bijective map $\phi$ from $\JJ_1 \rightarrow \JJ_2$  is said to be uniformly proper  if there exists a function $f: \natls \rightarrow \natls$ such that 
\begin{enumerate}
\item $d_{G_1} (J(L_i^1) ,J(L_j^1)) \leq n \Rightarrow d_{G_2} (\phi(J(L_i^1)) ,\phi(J(L_j^1))) \leq f(n)$ \\
\item  $ d_{G_2} (\phi(J(L_i^1)) ,\phi(J(L_j^1)))\leq n \Rightarrow d_{G_1} (J(L_i^1) ,J(L_j^1)) \leq f(n)$. 
\end{enumerate}
When $\JJ_i$ consists 
of all singleton subsets of $\Gamma_1, \Gamma_2$, we shall
refer to $\phi$ as a uniformly proper map from 
$\Gamma_1$ to $ \Gamma_2$.
\label{uproper}
\end{defn}

Our first main Theorem (combining Theorems \ref{main}, \ref{3mfld},
\ref{codim1}, \ref{even}, \ref{rk1} in the paper) is:\\

\begin{theorem} \label{omni} Let $n \geq 3$. Suppose 
$\JJ_i$ (for $i = 1,2$) are symmetric patterns of 
closed convex (or quasiconvex) sets in hyperbolic
space $\Hyp^{n} = \Hyp$, or more generally
uniform lattices in rank one symmetric spaces
of dimension $n$, $n \geq 3$. For $i= 1,2$, let $G_i$ be the
corresponding  cocompact group
of isometries.
Let $H_i \subset G_i$ be an 
 infinite index quasiconvex subgroup stabilizing the limit set
of some element of $\JJ_i$ and 
satisfying {\bf one} of the following
conditions: 

\begin{enumerate}
\item limit set of $H_i$ is a codimension
one topological sphere.
\item  limit set of $H_i$ is an even dimensional topological sphere.
\item  $H_i$ is a codimension one {\it duality group}. This generalizes 
(1). In particular, if  $n = 3$, $H_i$ could be {\it any} freely 
indecomposable subgroup of $G_i$. 
\item  $H_i$ is an odd-dimensional  Poincar\'e Duality group $PD(2k+1)$.
This generalizes (2). 

\end{enumerate}

Then any uniformly proper bijection between
$\JJ_1$ and $\JJ_2$ is induced by a hyperbolic isometry.
\end{theorem}

To prove Cases (1) and (2) we shall use the classical Brouwer 
and Lefschetz fixed point theorems respectively. To generalize these 
to Cases (3) and (4) we shall use tools from  the algebraic
topology  of generalized (or homological) manifolds.

Next suppose 
$\JJ$  is a  symmetric pattern of 
closed convex sets in  $ \Hyp$ as in Theorem \ref{omni}. For convenience
suppose that elements of $\JJ$ are $\epsilon$-neighborhoods of
convex hulls of limit sets of elements of $\JJ$, so that they are
strictly convex and $G$-equivariant.
Let $\phi$ be a uniformly proper bijection
from $\JJ$ to itself. Then Theorem \ref{omni} shows that
$\phi$ is induced by a hyperbolic isometry $f$. Consider the pattern
of geodesic segments perpendicular to elements of $\JJ$ at their end-points.
This collection is invariant under $G$ and there can only be a finite
number of such segments of bounded total length
 inside any bounded ball. Hence the subgroup of isometries
of $Isom ( \Hyp )$ preserving this pattern is discrete and contains
$G$ as a finite index subgroup. This proves the following
 Corollary  that is by now (after \cite{schwarz-inv} )
a standard consequence
of such pattern rigidity statements as Theorem \ref{omni}.

\begin{cor}
Suppose 
$\JJ$  is a  symmetric pattern of 
closed convex sets in hyperbolic
space $\Hyp^{n} = \Hyp$ or more generally
uniform lattices in rank one symmetric spaces
of dimension $n$, $n \geq 3$, as in Theorem \ref{omni} and $G$ the
associated cocompact group of isometries. Then
the subgroup of the quasi-isometry group
$QI( \Hyp )$ that coarsely preserves $\JJ$ contains $G$ as a subgroup
of finite index.
\label{commen}
\end{cor}

 More generally,
the pattern rigidity Theorem \ref{omni}
 goes through   for quasiconvex subgroups with
disconnected limit sets, at least one of whose components
has a stabilizer $H^{\prime}$ of the form (1), (2), (3)
 or (4) in Theorem \ref{omni}
above. For details, see Corollary \ref{disconnected} in this paper.
Theorem \ref{omni} combined with Corollary
\ref{disconnected} implies further that 
{\it pattern rigidity holds for all quasi
convex subgroups of hyperbolic 3-manifolds that are not virtually
free}. 

Similar extensions hold for quasiconvex subgroups
$H$ when some finite intersection of conjugates
$\bigcap_{i=1\cdots k} g_i H g_i^{-1}$ is of the form
(1), (2), (3)
 or (4) in Theorem \ref{omni}
above. For details see Corollary \ref{intersection} 
in this paper.

Combining this with the main theorem
of Mosher-Sageev-Whyte \cite{msw2} (to which we refer
for the terminology) we get the following QI-rigidity
theorem.

\smallskip

{\bf Theorem \ref{qirig} }
{\it Let $\GG$ be a finite, irreducible graph of groups with associated
Bass-Serre tree $T$ of spaces  
such that no depth zero raft of $T$ is a line.
Further suppose that the vertex groups are fundamental
groups of compact hyperbolic $n$-manifolds, or more generally
uniform lattices in rank one symmetric spaces
of dimension $n$, $n \geq 3$, and 
 edge groups are all of exactly one the following types:\\
a) A  {\it duality group} of codimension one
in the adjacent vertex groups.  In this case we require in addition 
that the crossing graph condition  of Theorems 1.5, 1.6
of  \cite{msw2} be satisfied
and that $\GG$ is of finite depth.\\
b) An odd-dimensional  Poincar\'e Duality group $PD(2k+1)$ with $2k+1 \leq n-1$.\\

\smallskip

 If $H$ is a finitely generated group quasi-isometric to $G = \pi_1 ( \GG )$ then $H$ splits as
a graph $\GG^{\prime}$ of groups whose depth zero vertex groups are commensurable to those of $\GG$ and
whose edge groups and positive depth vertex groups are respectively quasi-isometric to
groups of type (a),  (b).}

\subsection{Outline and Sketch} \label{outline}

\noindent {\bf Outline:} In Section 2, we describe some general
properties of limit sets of quasiconvex subgroups of hyperbolic
groups and recall some theorems from \cite{mahan-relrig}.
In Section 3, we recall some of the foundational work of Schwartz
from \cite{schwarz-inv} and describe some generalizations
that we shall use in this paper. Section 4 is the heart of the paper.
We reduce the problem of pattern rigidity to finding fixed points of
certain maps, and then proceed to apply classical fixed point theorems
(Brouwer and Lefschetz) to limit sets that are either spheres
of codimension one, or of even dimension. We generalize these results
to quasiconvex
Duality subgroups of dimension $n-1$ and quasiconvex $PD(2k+1)$
subgroups. For this we need some tools from the algebraic topology
of homology manifolds. In Section 5, we describe further generalizations
of these results to quasiconvex subgroups with disconnected
limit sets as well as subgroups with certain intersection properties.
We also combine these results with the main Theorem
of \cite{msw2} by Mosher-Sageev-Whyte to obtain QI-rigidity results.

\smallskip

\noindent {\bf Sketch of Proof:} We describe in brief the various steps
involved in the proof. \\
1) Uniformly proper pairings come from
quasi-isometries \cite{mahan-relrig}.\\
2) Use Mostow-Sullivan-Schwartz zooming in (cf Lemma \ref{zoom})
at a point of differentiability 
and non-conformality
to get an `eccentric' map $A$ on the boundary pairing limit sets.
The `eccentric' map is obtained by  pre- and
post-composing a linear map of Euclidean space
(thought of as a sphere minus the North pole)
with conformal maps of the sphere. This `zoom-in,
zoom-out' step is really quite classical and goes back to
Mostow \cite{mostow-pihes}. This was refined
by Sullivan \cite{sullivan-rigidity} and adapted to the present
context by Schwartz \cite{schwartz_pihes} \cite{schwarz-inv}. \\
3) Fix a particular limit set which is taken to another fixed limit
set under the
pairing. Zoom in using the stabilizer of the first limit set, act by $A$
and zoom out using the stabilizer of the second limit set. This step
differs from the corresponding step in \cite{schwarz-inv} as follows.
Though we can by the Generalized Eccentricity Lemma \ref{eccentric},
zoom in using powers of the same element, we cannot necessarily
zoom out using powers of the same element
 (as in \cite{schwarz-inv}). This makes the step technically more complicated
and we use a {\it generalized zoom-in zoom-out  Lemma} (Lemma \ref{contdiff}) 
to address this difficulty. \\
4) Get a sequence of rational functions that leave invariant a finite collection
of limit sets.\\
5) Apply Brouwer's Fixed Point Theorem (in the
codimension one sphere case) to get fixed points
in the ball
 bounded by the sphere;
and the Lefschetz Fixed Point Theorem (in the
even dimensional sphere limit set case) to get a fixed point
on the sphere limit set  itself. \\
6) Use some generalizations of Lefschetz
fixed point   theorem going back to work
of Lefschetz (himself), Felix Browder,
R. Thompson, R. Knill, R. Wilder along with a theorem
of Bestvina and  Bestvina-Mess to generalize Step 5 to
Duality and Poincar\'e duality groups. 

\subsection{Acknowledgments} The authors would like to thank
Marc Bourdon and Herve Pajot for beautiful lectures on
quasiconformal geometry and geometric function theory. We would
also like to thank Mladen Bestvina and Kalyan Mukherjea for  helpful
correspondence and discussions. We also thank the referee for several helpful
comments.

\section{Preliminaries}

\begin{rmk} {\rm 
A folklore fact that we shall be using is that for discrete subgroups $G$ of isometries of
real hyperbolic space ${\mathbb{H}}^n$, convex cocompactness is equivalent to quasiconvexity.
The forward implication is a consequence of the Milnor-Svarc Lemma (cf. \cite{gromov-hypgps}).

One way to prove the converse implication is to use the fact that ${\mathbb{H}}^n$
is projectively flat. Hence the convex hull of a finite set of
points is the union of the convex hulls of its $n+1$-tuples. Now let $\Lambda$ be the limit set of $G$ and let $CH_0$
denote the union of all ideal geodesic $n+1$-simplices with vertices in $\Lambda$. Then $CH_0$ is dense in the (closed)
convex hull $CH$ of $\Lambda$. Now, if an orbit $G.x$ is quasiconvex in ${\mathbb{H}}^n$, then the inclusion of $G.x$
into ${\mathbb{H}}^n$ extends continuously to a homeomorphism from the (Gromov) hyperbolic
boundary $\partial G$ to $\Lambda$. Hence, by $\delta$-hyperbolicity of $G$ and quasiconvexity, it follows
that for any $\{ x_0 , \cdots , x_n \} \subset \Lambda$ and  any $p$ in the ideal simplex with
vertices $\{ x_0 , \cdots , x_n \}$, $d(p, G.x)$ is uniformly bounded, where $d$ is the usual
metric on  ${\mathbb{H}}^n$. It follows that $CH_0$ and hence $CH$ lies in a bounded Hausdorff neighborhood of $G.x$.
Therefore the quotient $CH/G$ is compact. } \label{cc} \end{rmk}

\smallskip

\noindent{\bf Limit Sets and Pairings}\\
Let $G$ be a hyperbolic group. $\partial G$ will denote its boundary equipped with a visual metric.
Any fixed point of a hyperbolic element on $\partial G$ is called
a pole.
$\partial^2 G$ will denote the set of unordered pairs of distinct points on $\partial G$ with the topology
inherited from $\partial G$. A {\it pole-pair} is a pair of points $(x,y) \in \partial^2 G$ corresponding to
the fixed points of a hyperbolic element of $G$.

\begin{lemma}  (Gromov \cite{gromov-hypgps}, 8.2G, p.213) Pole-pairs are dense in $\partial^2 G$ and more generally, if  $H$ be a  finitely generated
group of isometries acting on hyperbolic space $\Hyp$ with limit
set $\Lambda$ then pole-pairs are dense in  $\partial^2 \Lambda$.
\label{poles}
\end{lemma}

The next Lemma is a consequence of the fact that the
action of a finitely generated
group of isometries of a hyperbolic metric space $\Hyp$ 
acting on the limit set is a convergence group action.

\begin{lemma} Let $H$ be a  finitely generated
group of isometries acting on hyperbolic space $\Hyp$ with limit
set $\Lambda$. Then for all
$(x,y) \in \partial^2 \Lambda$, there exists a sequence of hyperbolic isometries $T_i \in H$ with attracting
(resp. repelling) fixed points $x_i$ (resp. $y_i$) such that
$x_i \rightarrow x$,  $y_i \rightarrow y$, and
 the translation length of $T_i$ tends to
$\infty$.
\label{attract}
\end{lemma}

\noindent {\bf Proof:} We choose pole-pairs $(x_i,y_i)$ converging to $(x,y) \in \partial^2 \Lambda$ by Lemma \ref{poles}.
Let $T_i$ be a hyperbolic isometry
in $G$ with attracting fixed point $x_i$ and repelling fixed point $y_i$.
Choosing appropriately large powers $T_i^{n_i}$ of $T_i$, we are through. 
$\Box$ 

Since the orbit of an open subset of $\partial G$ under $G$ is the whole of $\partial G$, it follows that the limit
set $L_H$ of 
any infinite index quasiconvex subgroup $H$ of $G$ is nowhere dense in $\partial G$. Assume for simplicity
that $H = Stab (L_H)$. Then for all $g \in G \setminus H$, $gHg^{-1} \cap H$ is an infinite
index quasiconvex subgroup of $H$ (by a Theorem of Short \cite{short}) and hence its limit
set is nowhere dense in $L_H$. As $g$ ranges over $g \in G \setminus H$, we get a countable collection
of nowhere dense subsets of $L_H$. The next Lemma follows.

\begin{lemma} \label{baire} Suppose that $H = Stab (L_H)$ is a quasiconvex subgroup of
a hyperbolic group $G$. 
For all $p \in L_H$ and all $\epsilon > 0$ there exists $x \in L_H$ such that $d(p,x) < \epsilon$ and $L_H$ is the
unique translate of $L_H$ to which $x$ belongs, i.e. if $ x \in L_H \cap gL_H$ then $g \in  H$.
\end{lemma}

\begin{lemma} Suppose that $H = Stab (L_H)$ is a quasiconvex subgroup of
a hyperbolic group $G$. Let $U \subset \partial G$ be an open subset and let $\epsilon > 0$. Then there exists
a finite collection of points $x_1, \cdots , x_n \in U$ such that \\
\begin{enumerate}
\item  $\{   x_1, \cdots , x_n  \}$ is an $\epsilon$-net in $U$ 
\item  $x_i \in L_i = g_i L_H$ for some $g_i  \in G \setminus H$ and $L_i$ is the
unique translate of $L_H$ to which $x$ belongs.
\end{enumerate}
\label{net}
\end{lemma}

\noindent
{\bf Proof:}
This follows from Lemma \ref{baire} and the fact that the union of all the translates of $L_H$ under $G$ is 
dense in $\partial G$. $\Box$

\begin{defn} A point that belongs to a unique translate of $L_H$ will
be called a {\bf unique point}.
\end{defn}

In \cite{mahan-relrig} we showed the following:

\begin{theorem} \label{qipairs} (Theorem 3.5 of \cite{mahan-relrig})\\
Let $\phi$ be a uniformly proper (bijective, by definition) map from $\JJ_1 \rightarrow \JJ_2$. There exists a quasi-isometry $q$ from $\Gamma_1$ to $\Gamma_2$ which pairs the sets $\JJ_1$ and $\JJ_2$ as $\phi$ does.
\end{theorem}

\begin{prop} {\bf Characterization of Quasiconvexity} ( Prop 2.3 of \cite{mahan-relrig} )\\
Let $H$ be a subgroup of a hyperbolic group $G$ with limit set $\Lambda$. Let $\LL$ be the collection of translates of $\Lambda$ (counted with multiplicity)
 by elements of distinct cosets of $H$ (one for each coset). Then $H$ is quasiconvex if and only if $\LL$ is a discrete subset of 
$C_c^0(\partial G)$,
where $C_c^0(\partial G)$ denotes the collection of compact subsets of
$\partial G$ with more than one point equipped with the Hausdorff metric.
\label{qc=discrete}
\end{prop}

Finally, combining Lemmas \ref{poles} and \ref{baire} along with Proposition
\ref{qc=discrete}, we get

\begin{cor} {\bf Generic Pole-pairs}
Suppose that $H = Stab (L_H)$ is a quasiconvex subgroup of
a hyperbolic group $G$. Identify $L_H$ with the boundary $\partial H$ of
$H$.
For all $(p,q) \in \partial^2 H$ and all $\epsilon > 0$ there exists a pole-pair
$(x,y) \in \partial^2 H$ such that $d(p,x) < \epsilon$,
$d(q,y) < \epsilon$
 and $L_H$ is the
unique translate of $L_H$ to which $x$ (or $y$)
belongs, i.e. if $ x$ 
$ ( \rm{ or })$ $ y  \in L_H \cap gL_H$ then $g \in  H$.
\label{genericpole}
\end{cor}

\medskip

\noindent {\bf Small Homotopies}\\
We shall have need for the following fact \cite{bott-tu}.

\begin{lemma} \label{homotopy} Given a closed Riemannian
manifold $(M,d)$, there exists $\epsilon_1 > 0$ such that the following holds: \\
If $f$ is a self-homeomorphism
 such that $d(f(x), x) < \epsilon_1 $ for all $x \in M$, then $f$ is homotopic to the identity. 
\end{lemma}

\noindent {\bf Sketch of proof for smooth maps:}
Since $M$ is a compact Riemannian manifold, there is an $\epsilon_0>$
such that (modulo the natural identification of the normal bundle of
the diagonal, $D_M\subset (M\times M)$, with the tangent bundle of
$M$) tangent  vectors of length $\epsilon_0$ map  via the
exponential map diffeomorphically onto an open
tubular neighborhood  of $D_M$.
If $f: M\to M$ is sufficiently near to the identity
$1_M$ in the compact-open
topology, then the image of the graph of $f$ will lie in this tubular
neighborhood. The inverse of the exponential map identifies the graph of
$F$ with a section of the tangent bundle of $M$ and since any section
is homotopic to the zero section (which corresponds to the graph of
the identity map) we are through. $\Box$

\medskip

\noindent {\bf Boundaries}\\
We shall have need for the following Theorem of Bestvina and Mess \cite{bestvina-mess}.

\begin{theorem} \label{bestvina0} Boundaries $\partial G$
of PD(n) hyperbolic
groups $G$ are locally connected
homological manifolds (over the integers) with the homology
of a sphere of dimension $(n-1)$. Further, if $G$ 
acts freely, properly, cocompactly on a contractible complex $X$, then the natural compactification
$X \cup \partial G$ is an Absolute Retract (AR). \end{theorem}

\begin{rmk} {\rm In particular, $X \cup \partial G$ has the fixed point property, i.e. any continuous
map from $X \cup \partial G$ to itself has a fixed point.} \label{bestvina1}  \end{rmk}

\section{Differentiability Principles and Eccentric Maps}

\subsection{Differentiability Principles}
Let $ {\Hyp}^{n+1} = \Hyp$ denote hyperbolic ($n+1$)-space and 
$\bdy {\Hyp}^{n+1} = S^n_\infty$ denote the boundary sphere at infinity 
with the standard conformal structure (preserved by isometries
of $ {\Hyp}^{n+1} $). Let ${\E}^n = \E$ denote the Euclidean space 
obtained from $S^n_\infty$ by removing the point at infinity.

We recall  a certain {\it Differentiability Principle} 
from Schwartz's
paper \cite{schwarz-inv}. Suppose $h \co S^n_\infty \rightarrow
S^n_\infty$ is a homeomorphism fixing $0, \infty$ such that $dh(0)$
exists. Let $T_1, T_2$ be two contracting similarities 
(with possible rotational
components) of
$\E$ both fixing $0$.
For each pair $k_1, k_2$ of positive integers, Schwartz defines
the map $$h[k_1, k_2] = T_2^{-k_2} \circ h \circ T_1^{k_1}$$
and shows

\begin{lemma} (Lemma 5.3 of \cite{schwarz-inv} )\label{zoom}
Suppose that $K_1, K_2 \subset \E$ are compact subsets. Suppose that
$(k_{11}, k_{21}), (k_{12}, k_{22}), (k_{13}, k_{23}), \cdots $ is 
a sequence
of pairs such that \\
\begin{enumerate}
\item  $k_{1n} \rightarrow \infty$. 
\item  $h[k_{1n}, k_{2n}] (K_1) \cap K_2 \neq \emptyset$.
\end{enumerate}
Then on some subsequence $h[k_{1n}, k_{2n}]$ converges, uniformly on
compact sets, to a linear map.
\end{lemma}

\begin{rmk} {\rm Lemma \ref{zoom} can be slightly generalized by replacing the maps $T_1^{k_1}$ and 
$T_2^{k_2}$ by maps $T_{1k_1}$ and $T_{2k_2}$ such that the translation lengths of $T_{1k_1}$ and $T_{2k_2}$
tend to infinity as $k_1, k_2 \rightarrow \infty$. This is all Schwartz uses in his proof.}
\label{zoom-rmk} \end{rmk}

We shall need a generalization and weakening
of this to continuously differentiable functions.

\begin{lemma} \label{contdiff} {\bf Generalized Zoom-in Zoom-out} \\
Suppose $h \co S^n_\infty \rightarrow
S^n_\infty$ is continuously differentiable.
Let $T_{1n}, T_{2n}$ be sequences of hyperbolic Mobius transformations
such that their fixed point sets $\{ x_{1n}, y_{1n} \}$,
$\{ x_{2n}, y_{2n} \}$ satisfy 
\begin{enumerate}
\item  $x_{in} $  is the attracting fixed point of $T_{in}$, for $i= 1,2$.
\item  $y_{in} $  is the repelling fixed point of $T_{in}$, for $i= 1,2$.
\end{enumerate}
We further assume that $x_{1n} = 0$ and  $y_{1n} = \infty$ for all $n$. \\
\smallskip

Let 
$$h_n = T_{2n}^{-1} \circ h \circ T_{1n}$$
Suppose that
$(T_{11}, T_{21}), (T_{12}, T_{22}), (T_{13}, T_{23}), \cdots $ is a sequence
of pairs such that \\
a) The translation lengths of $T_{1n} \rightarrow \infty$. \\
b) There exists $\epsilon > 0$ such that
 $${\rm inf}_n ({\rm min} \{ d(h_n (0), h_n (1)),
d(h_n ( \infty ), h_n (1)), d(h_n (0), h_n ( \infty )) \} )\geq \epsilon$$
Then on some subsequence $h_n$ converges, uniformly on
compact sets, to a linear map post-composed with a conformal map.
\end{lemma}

\noindent {\bf Proof:} The sequence
$g_n = T_{1n}^{-1} \circ h \circ T_{1n}$ converges (up to sub-sequencing)
to a linear map by Lemma \ref{zoom} and Remark \ref{zoom-rmk}. 
Condition (b) in the hypothesis
guarantees that the ratio of
the translation lengths of $T_{1n}$ and $T_{2n}$ is bounded away from
both $0$ and $\infty$. Hence, (extracting a further subsequence
if necessary) $T_{2n}^{-1}  \circ T_{1n}$ converges to a
conformal map. The result follows. $\Box$

\subsection{Eccentric Maps}

\begin{defn} \cite{schwarz-inv}
Let $T$ be a real linear map of the Euclidean space ${\E}^n = \E$. Let 
$g_i$ (for $i=1,2$) be two conformal maps of $\bdy {\Hyp}^{n+1} = S^n_\infty$.
The map $\mu = g_2 \circ T \circ g_1^{-1}$ is said to be
an {\bf eccentric}
map if 
\begin{enumerate}
\item  $\mu$ preserves $E$ and fixes $0$.
\item  $\mu$ is differentiable at $0$. 
\item  $\mu$ is \underline{not} a real linear map.
\end{enumerate}
\end{defn}

Then  the {\it Eccentricity Lemma} (Lemma 2.2) of Schwartz 
\cite{schwarz-inv} generalizes to

\begin {lemma} {\bf Generalized Eccentricity Lemma} Let
$G_1, G_2$ be two groups acting freely, properly discontinuously
by isometries and cocompactly on $ {\Hyp}^{n+1} = \Hyp$.
Let  $H_1^0, H_2^0$ be quasiconvex subgroups of $G_1, G_2$ respectively
with limit sets $\Lambda_1^0, \Lambda_2^0$.
Let $\JJ_i^0$ (for $i = 1, 2$) be the set of translates of
joins (or convex hulls) 
of $\Lambda_i^0$. Let $q^0$ be a quasi-isometry pairing
$\JJ_1^0$ with $\JJ_2^0$. Assume that $h^0 = \partial q^0$ is not conformal.
Then there exist symmetric pattern of joins (or convex hulls) $\JJ_i$ (of limit sets 
$\Lambda_i$ abstractly homeomorphic to $\Lambda_1^0, \Lambda_2^0$)
and a quasi-isometry $q \co \Hyp \rightarrow \Hyp$ such that 
\begin{enumerate}
\item $q$ pairs the elements of $\JJ_1$ with those of $\JJ_2$ 
\item  $\mu = \partial q$ is an eccentric map. 
\item  The geodesic $\gamma = \overline{0 \infty}$ is a subset of some
$J_i \in \JJ_i$ for $i = 1, 2$. Further, the (translates
of the) limit sets $\Lambda_1$ and $\Lambda_2$ in which the end-points
${0,  \infty}$ lie is unique.  
\item ${0,  \infty}$ are poles for the action of the stabilizer
Stab ($\Lambda_1$) on $\Lambda_1$. 
\end{enumerate}
\label{eccentric}
\end{lemma}

\noindent {\bf Proof:}
The difference with Lemma 2.2 of
\cite{schwarz-inv} is in Conditions (3) and (4) above. Schwartz' proof 
proceeds by zooming in at a point of differentiability 
and non-conformality (taken to be the
origin)
of the quasiconformal map
$h^0$ to  obtain a linear map $h^{\prime}$
from $\E$ to itself in the limit. 
The sequence of maps used in zooming in come by
conjugating $h^0$ by $D^n$ where $D$ is a dilatation map with ${0,  \infty}$ as fixed points. Further $h^{\prime}$ is the boundary value 
of some quasi-isometry $q^{\prime}$ which pairs some symmetric pattern 
of joins $\JJ_1^{\prime}$ with $\JJ_2^{\prime}$.
This step
goes through verbatim. 

Next,
by Lemma \ref{baire} 
there exist pairs of points $\alpha , \beta$ 
on some limit set $\Lambda_1$ of an element of $\JJ_1^{\prime}$
such that $\Lambda_1$ is unique, i.e.  $\alpha , \beta$ do not belong
to any other limit set $\Lambda_1^{\prime}$ 
of an element of $\JJ_1^{\prime}$. Also,
by Corollary \ref{genericpole} the pair
$( \alpha , \beta )$ can be taken as
  pole-pair for the action of
the stabilizer
Stab ($\Lambda_1$) on $\Lambda_1$.  Since 
$q^{\prime}$ pairs the symmetric pattern 
of joins $\JJ_1^{\prime}$ with $\JJ_2^{\prime}$, 
$h^{\prime}( \alpha ) = \alpha^{\prime} , h^{\prime}(\beta )
= \beta^{\prime}$
belong to some unique $\Lambda_2$, i.e. 
$h^{\prime}( \alpha ), h^{\prime}(\beta )$ 
do not belong
to any other limit set $\Lambda_2^{\prime}$ 
of an element of $\JJ_2^{\prime}$. Let $g_j$ be chosen in such a 
way that $g_1$ (resp. $g_2$) maps ${0,  \infty}$ to $ \alpha , \beta$
(resp. $\alpha^{\prime}, \beta^{\prime}$ ) respectively. 

Then $\mu = g_2 \circ h^{\prime} \circ g_1^{-1}$
and $q = g_2 \circ q^{\prime} \circ g_1^{-1}$
are the required maps. $\Box$

We shall need the following `Zariski-density'
property of eccentric maps due to Schwartz
\cite{schwarz-inv}.

\begin{lemma} (Corollary 5.2 of \cite{schwarz-inv}). Let $U \subset
\E$ be an open subset. Then there is a constant $\delta
= \delta (U) > 0$ such that if two eccentric maps agree on a 
$\delta$-dense subset of $U$, then they agree everywhere.
\label{zariski}
\end{lemma}

\section{Pattern Rigidity}

\subsection{Scattering}

For $i = 1, 2 $, 
let $F_i$ be a (compact) fundamental domain for the action 
of $H_i = Stab(L_{H_i})$ on the domain of discontinuity $\Omega_{i}$
of $H_i$ (see Remark \ref{cc}). Let $Q_i$ be the quotient of $\Omega_i$ by $H_i$. Let $\Pi_i
\co \Omega_{i} \rightarrow Q_i$ 
be the covering map. Recall, by Lemma \ref{eccentric}
 that $0, \infty$ form a 
 pole-pair for the action of $H_1$ on $L_{H_1}$. Next, suppose that
 we have \\
\begin{enumerate}
\item An eccentric map $\mu$ \\
\item A subset $\Sigma \subset Q_1$ \\
\item A neighborhood $S \subset E$ of $0$ \\
\end{enumerate}

Define $\Psi (\mu , \Sigma , S) = \Pi_2 \circ \mu (S \cap \Pi_1^{-1} ( \Sigma ) ) \subset Q_2$. 

\begin{lemma} {\bf Scattering Lemma:} Independent of $\mu$ there is a constant $\delta_0 > 0$ such that if $S \subset E$ is  any
neighborhood of $0$,  and $\Sigma \subset Q_1$ is $\delta_0$-dense, then 
$\Psi (\mu , \Sigma , S)  \subset Q_2$ is an infinite set. 
\label{scatter}
\end{lemma}

\noindent {\bf Proof:} Though we shall follow the broad scheme
of the proof of Lemma 2.3 of Schwartz \cite{schwarz-inv},  technically our proof
 will be quite
a bit more involved as we shall first 
use the Generalized Zoom-in Zoom-out Lemma \ref{contdiff} and then Lemma 
\ref{zoom} (and not Lemma \ref{zoom} directly as in \cite{schwarz-inv}). 
In particular, steps (1) and (2) below will be different, while
step (3) will be the same as in \cite{schwarz-inv}.

Let $\Sigma_0 = \Pi_1^{-1} (\Sigma ) \cap F_1$. Let $S$ be an open neighborhood of $0$. 

There exists
by the Generalized Eccentricity
Lemma \ref{eccentric} a sequence of hyperbolic Mobius transformations
$T_{1n} \in H_1$ 
such  that the fixed point sets $\{ x_{1n}, y_{1n} \}$,
 satisfy 
\begin{enumerate}
\item $x_{1n} = 0$ is the attracting fixed point of $T_{1n}$. 
\item $y_{1n} = \infty$ is the repelling fixed point of $T_{1n}$. 
\item The hyperbolic isometries corresponding to 
$T_{1n}$, form an unbounded set in $PSL_2(C)$.
\item $T_{1n} (F_1) \subset S$.
\end{enumerate}

Condition (4) follows from (1) and substituting $T_{1n}$ by
 large enough powers of $T_{1n}$
if necessary.

\smallskip

\noindent {\bf Step 1: Choosing $T_{2n} \in H_2$}\\
The first step is to choose $T_{2n} \in H_2$.
 $T_{1n} (\Sigma_0 ) \subset \Pi_1^{-1} (\Sigma ) \cap S$.
First choose a point $w \in F_1$. Let $P(w)$ denote the foot of the perpendicular from $w$ to $(0, \infty ) = (x_{1n}, y_{1n})$. Then
$T_{1n} (w) \in S$. Map the tripod with vertices $(0, \infty , T_{1n} (w) )$
over by $\mu$. Recall that $\mu (0) = 0, \mu (\infty ) = \infty$.
Then $(0, \infty , \mu \circ T_{1n} (w) )$ form the vertices
of a {\em uniform}
$K$-quasitripod in $\Hyp$ with centroid $\mu \circ T_{1n} (P(w))$.
Choose $T_{2n} \in H_2$ such that  $T_{2n}^{-1} \circ \mu 
\circ T_{1n} (P( w) ) $ lies in a fixed fundamental domain for the action 
of $H_2$ on the convex hull $CH( L_{H_2} )$ of the limit set
$ L_{H_2}$. (We could equally well have chosen the join 
$J( L_{H_2} )$ of the limit set
$ L_{H_2}$.) Then, automatically, the three points
$T_{2n}^{-1} (0)$, $T_{2n}^{-1} ( \infty )$, 
$T_{2n}^{-1} \circ \mu 
\circ T_{1n} ( w)  $ satisfy the hypotheses of Lemma \ref{contdiff},
i.e. there exists $\epsilon > 0$ such that 
the three points
$T_{2n}^{-1} (0)$, $T_{2n}^{-1} ( \infty )$, 
$T_{2n}^{-1} \circ \mu 
\circ T_{1n} ( w)  $ are at a distance of at least $\epsilon$
from each other on the sphere (uniformly for all $n$). 
By Lemma \ref{contdiff}, (up to extracting a subsequence)
$T_{2n}^{-1} \circ \mu 
\circ T_{1n}$ converges to a map $\psi \circ L$, where $\psi$ is conformal
and $L$ is linear. 

\smallskip

\noindent {\bf Step 2: The sequence $T_{2n}^{-1} \circ \mu 
\circ T_{1n}$ consists of infinitely many distinct elements}\\
We would like to conclude that there are {\it infinitely many
distinct maps} in the sequence $T_{2n}^{-1} \circ \mu 
\circ T_{1n}$ converging to a map $\psi \circ L$. Suppose not.
Then the sequence
of maps 
$T_{2n}^{-1} \circ \mu 
\circ T_{1n}$ is eventually constant and equal to $\psi \circ L$. In particular, since $\mu, T_{1n}, L $ all fix $0, \infty$, it follows
that $T_{2n}^{-1} (0) = \psi (0)$ and $T_{2n}^{-1} ( \infty ) = \psi (
 \infty )$ for all $n$. But then 
\begin{enumerate}
\item  $T_{2m} \circ T_{2n}^{-1} (0) = 0$
and $T_{2m} \circ T_{2n}^{-1} ( \infty ) =  \infty $ for all $m, n$. 
\item   $T_{2n} \in H_2$ for all $n$.
\end{enumerate}

 Since the collection
of elements of the form $T_{2m} \circ T_{2n}^{-1}$ is
infinite, it follows that  the set of elements in $H_2$
fixing $0, \infty$ is virtually infinite cyclic. Thus, $(0, \infty )$ form
a  pole-pair for the action of $H_2$ on $L_{H_2}$. Let $\CC$
be an infinite cyclic subgroup of $H_2$ fixing $0, \infty$.
In this case,
we modify the sequence $T_{2n}$ by choosing these to be elements
of  $\CC$ satisfying the hypotheses of Lemma \ref{zoom}. Then 
$\mu_n = T_{2n}^{-1} \circ \mu 
\circ T_{1n}$ converges to a linear map by Lemma \ref{zoom}.
Since $\mu$ is eccentric, so is $\mu_n$ and
we may assume that $\mu_n \rightarrow \mu^{\prime}$, a linear map.
Hence in either case, we can conclude that 
the sequence $\mu_n = T_{2n}^{-1} \circ \mu 
\circ T_{1n}$ of maps consists of infinitely many distinct elements 
converging either to a map of the form $\psi \circ L$ (with $\psi$ conformal and $L$ linear) 
or simply a linear map $L$.

\smallskip

\noindent {\bf Step 3: Using  Zariski density  }\\
The rest of the proof follows that of \cite{schwarz-inv}. Define
$V = \bigcup_{n=1}^{\infty} \mu_n (\Sigma_0 )$. 

\begin{claim}
$V$ is  infinite and  $\overline{V} \subset \Omega_2$.
\end{claim}

\noindent {\bf Proof of Claim:}  Since $\mu_n \rightarrow \mu$,
$V$ is bounded and  $\overline{V} \subset \Omega_2$. In particular,  $\overline{V}$
is contained in the union of finitely many translates of the compact fundamental domain $F_2$. By Step (2)
  there are infinitely many distinct maps in the sequence. If $V$ is finite,
then only finitely many choices are there for $\mu_n (\Sigma_0 )$ and hence
by Lemma \ref{zariski}, there are only finitely many choices
for $\mu_n$. This contradiction proves that $V$ is infinite. $\Box$

Since $V$ is  infinite and $\overline{V}$
is contained in the union of finitely many translates of the compact fundamental domain $F_2$, 
it follows that $\Pi_2 (V)$ is
infinite. But $\Pi_2 (V) \subset \Psi (\mu , \Sigma , S)  \subset Q_2$.
Hence $\Psi (\mu , \Sigma , S)$ is infinite. $\Box$

\subsection{Pattern Rigidity: Topological Spheres} \label{top}
In this subsection we shall prove pattern rigidity for symmetric patterns of 
closed convex (or quasiconvex) sets in hyperbolic
space $\Hyp^{n+1} = \Hyp$ such that the limit sets are topological spheres
(of either codimension one or of even dimension). The techniques used are from fixed point theory.
In the next two subsections, we shall generalize this, 
to quasiconvex
subgroups of 3-manifolds with connected limit sets, to codimension one quasiconvex
duality subgroups, 
and  to closed limit sets whose stabilizers are $PD(2n+1)$
quasiconvex
 subgroups. 
These will be generalizations of Theorem \ref{main} (a) and Theorem \ref{main} (b)
respectively. 
The technicalities for these generalizations
are postponed  for ease of exposition.

Recall that a point that belongs to a unique translate of $L_H$ is
 called a {\bf unique point}.

\begin{theorem} \label{main}
Let $n \geq 2$. Let
$\JJ_i^0$ (for $i = 1,2$) be symmetric patterns of 
closed convex (or quasiconvex) sets in hyperbolic
space $\Hyp^{n+1} = \Hyp$ such that the limit sets
of $\JJ_i^0$ are either \\
a) topological spheres of dimension $(n-1)$, {\it OR} \\ 
b) even-dimensional topological spheres. \\
Then any proper bijection $\phi$ between
$\JJ_1^0$ and $\JJ_2^0$ is induced by a hyperbolic isometry.
\end{theorem}

\noindent {\bf Proof:} By Theorem \ref{qipairs},
there is a quasi-isometry $q^0$ that pairs the convex (or quasiconvex) sets 
$\JJ_i^0$ as $\phi$ does. 

Suppose that $h^0 = \partial q^0$ is not conformal.
Then by the {\bf Generalized Eccentricity Lemma \ref{eccentric} }
there exist, for $i = 1, 2$,
 symmetric patterns of convex, or quasiconvex sets
 $\JJ_i$ (with limit sets 
$\Lambda_1$ abstractly homeomorphic to $\Lambda_2$)
and a quasi-isometry $q \co \Hyp \rightarrow \Hyp$ such that 
\begin{enumerate}
\item $q$ pairs the elements of $\JJ_1$ with those of $\JJ_2$ \\
\item $\mu = \partial q$ is an eccentric map. \\
\item The geodesic $\gamma = \overline{0 \infty}$ is a subset of some
$J_i \in \JJ_i$ for $i = 1, 2$. Further, the (translates
of the) limit sets $\Lambda_1$ and $\Lambda_2$ in which the end-points
${0,  \infty}$ lie is unique.
\end{enumerate}

Let $\delta_0$ be as in Lemma \ref{zariski}.
Pick points as per Lemma \ref{net} to get a $\delta_0$-net $\Sigma$
in the
interior of the fundamental domain $F_1$ of the action of $H_1$
on its domain
of discontinuity $\Omega_1$, consisting of unique points. Let $S$ be an open neighborhood
of $0$. Then, by the Scattering Lemma \ref{scatter}
$\Psi (\mu , \Pi_1 (\Sigma ), S) = \Pi_2 \circ \mu (S \cap \Pi_1^{-1} (\Pi_1( \Sigma )) ) \subset Q_2$ is infinite.

However, since $\Sigma = \{ x_1, \cdots x_n \}$ 
is finite, and since its points belong to
unique limit sets, ($x_i \in L_i$ say) 
there is an upper bound on the distance
of $J(L_{H_1})$ from $J(L_i)$. Since $q$ is a quasi-isometry, there
is an upper bound on the distance of $J(L_{H_2})$ from
$\phi (J(L_i))$. Hence, modulo the action of $H_2$, there are 
only finitely many choices for $\phi (J(L_i))$.

Since $\Psi (\mu , \Pi_1 (\Sigma ), S) = \Pi_2 \circ \mu (S \cap \Pi_1^{-1} (\Pi_1( \Sigma )) ) \subset Q_2$ is infinite, it follows that
there exists (after subsequencing again) $T_{in} \in H_i$ for $i = 1,2$
such that 
\begin{enumerate}
\item If
$\mu_n = T_{2n}^{-1} \circ \mu \circ T_{1n}$, then $\mu_n (L_{H_1})
= L_{H_2}$
and for some $L_i = L_1$ (say, without loss of generality)
$\mu_n (L_{1})
= L_{2}$ is a fixed limit set. This follows from the fact that the $x_i$'s are unique points. Also note that we can arrange 
that the visual diameters of $L_1, L_2$ are smaller than any
pre-assigned $\epsilon_0$. 
\item The attracting (resp. repelling) fixed point $x_{1n}$ (resp. 
$y_{1n}$) are $0$ (resp. $\infty$).
\item  The hyperbolic isometries corresponding to 
$T_{1n}$ form an unbounded set in $PSL_2(C)$.
\item  $\mu_n$ restricted to ${L_1}$ are distinct maps as 
$\Psi (\mu , \Pi_1 (\Sigma ), S)$ is infinite. In particular,
$\mu_n$'s are distinct maps.
\item  $\mu_n \rightarrow \mu^\prime$, where $\mu^{\prime}$ is either
a real linear
map or a real linear
map post-composed with a conformal map  where 
the linear factor of the map $\mu^{\prime}$  is not a similarity, but
 continues to satisfy property (1). 
\end{enumerate}

Further, by Proposition \ref{qc=discrete}, if we fix any {\it finite}
 collection of translates, $L_{11}, \cdots , L_{1m}$, of the limit set
$L_{H_1}$, then the (ordered tuple) 
  $\mu_n (L_{11}), \cdots ,
\mu_n (L_{1m})$ is eventually constant. Hence
for $n,l$ sufficiently large, $\mu_l^{-1} \circ \mu_n$ maps
$L_{1j}$ to itself for $j = 1, \cdots , m$. 

The argument so far does not use any special topological
property of the limit sets.
We summarize our conclusions in the Remark below.

\begin{rmk} \label{invariantsets}  {\rm We have 
 shown that given an eccentric 
map pairing symmetric patterns $\JJ_i$
of convex (or quasiconvex) sets, there 
exists
\begin{enumerate}
\item  a sequence of eccentric maps $\mu_j \rightarrow \mu^{\prime}$ uniformly
on compact sets, where
$\mu^{\prime}$ is \\
a) either a linear map that is {\bf not} a similarity.\\
b) or a real linear
map post-composed with a conformal map  where 
the linear factor of the map   is not a similarity.\\
\item  $\mu_j$'s pair $\JJ_1$ with $\JJ_2$ \\
\item For any finite collection  $\LL$ of limit sets
of elements of $\JJ_1$, there exists a positive integer $N$, such that
$\mu_n (L) = \mu_l (L)$ for all $L \in \LL$ and $n, l \geq N$. \\
\item  $\mu_j$'s are {\bf distinct} eccentric maps.
\end{enumerate}
}
\end{rmk}

We now deal with the two cases of the Theorem separately.

\smallskip

{\bf Case a:} Limit sets
of $\JJ_i^0$ are topological spheres of dimension $(n-1)$. \\
Since each $L_{1i}$ is a topological
sphere of codimension one and $J_{1i}$ is a convex set
(for quasiconvex sets, we take the convex hull), the compactification of
a small $\epsilon$-neighborhood,
$N_\epsilon (J_{1i} )$ obtained by adjoining $L_{1i}$ is a strong deformation
retract of the whole compactified ball
$\D = \Hyp \cup S^n_\infty$. In particular, if $\Omega$ is one of the two components of the domain
of discontinuity of $Stab(L_{1i})$, then $\Omega \cup L_{1i} = D_{1i}$ is an AR by Theorem \ref{bestvina0}
and hence satisfies the fixed-point property (Remark \ref{bestvina1}).

For $n,l$ sufficiently large, $\mu_l^{-1} \circ \mu_n$ maps
$D_{1i}$ to itself for $i = 1 \cdots m$. By {\it Brouwer's fixed point
Theorem} (Remark \ref{bestvina1}), there exist $x_{1i} \in D_{1i}$, such that 
$\mu_n  (x_{1i}) =  \mu_l (x_{1i})$ for $i = 1 \cdots m$.
Now, by Remark \ref{invariantsets}
 above, we can choose $L_{1i}$ of sufficiently
small diameter such that for any $x_{1i} \in D_{1i}$, the collection
$\{$ $x_{11}, \cdots , x_{1m}$ $\}$ is an $\epsilon_0$-net
in $S^n$, where $\epsilon_0$ is as in Lemma \ref{zariski}.
Hence, by Lemma \ref{zariski}, $\mu_n   =  \mu_l $. This contradicts Condition (4) of Remark \ref{invariantsets}
above and proves Case (a) of the Theorem. $\Box$

\smallskip

{\bf Case b:} Limit sets
of $\JJ_i^0$ are topological spheres of even dimension. \\

By Lemma \ref{homotopy} and   Remark \ref{invariantsets} it follows that given
 any finite collection $\LL$ of limit sets, there exists a positive integer $N$ such that
for all $n, l \geq N$, and all ${L}_i \in \LL$, \\
\begin{enumerate}
\item $ \mu_l^{-1} \circ \mu_n (L_i )  = L_i$
\item $ \mu_l^{-1} \circ \mu_n$ restricted to $L_i$ is homotopic
to the identity. Hence, 
the Lefschetz number of $ \mu_l^{-1} \circ \mu_n$ restricted to $L_i$ is
equal to the Euler characteristic of $L_i$ 
\end{enumerate}

Since each $L_i \in \LL$ is an even-dimensional sphere, the Euler
characteristic of $L_i$ is $2$, in particular non-zero.
By the Lefschetz fixed point Theorem there exists $x_i \in L_i$ such that
$\mu_n (x_i) = \mu_l (x_i)$. 

The rest of the proof is as in Case (a) above. 
By Remark \ref{invariantsets}, we can choose $L_{i}$ of sufficiently
small diameter such that for any $x_{i} \in L_{i}$, the collection
$\{$$x_{1}, \cdots , x_{m}$$\}$ is an $\epsilon_0$-net
in $S^n$, where $\epsilon_0$ is as in Lemma \ref{zariski}.
Hence, by Lemma \ref{zariski}, $\mu_n   =  \mu_l $. This contradicts Condition (4) of Remark \ref{invariantsets}
above and proves Case (b) of the Theorem. $\Box$

\begin{rmk} {\rm The Proof of Case (b) when specialized to dimension
 zero (i.e. $S^0$ limit sets) is exactly the one given by Schwartz
in \cite{schwarz-inv}. To see this note that
  the existence of a fixed
point of a map from $S^0$ to itself that is `close to the identity'  
(hence equal to the identity) is clear. }
\end{rmk}

\subsection{3 manifolds and Codimension one Duality Subgroups}

\begin{rmk} \label{genlzn} {\rm In our proof of Theorem \ref{main},
we have  used the fact that the limit sets
are spheres in a mild way. In case (a) we used them to construct invariant 
Absolute Retracts bounded by these spheres. 
After this, the proof of both Case (a) and Case (b) end up
using the Lefschetz fixed point Theorem. We have 
used  the following  facts:
\begin{enumerate}
\item  Euler characteristic of each invariant limit 
set $L$ is non-zero. 
\item  A map that moves each point of $L$ through a small distance
is homotopic to the identity. 
\item  The Lefschetz fixed point Theorem holds for $L$. 
\end{enumerate}
}
\end{rmk}

We generalize Theorem \ref{main} (a) now to quasiconvex
subgroups of 3-manifolds with connected limit sets.

\begin{theorem}
 \label{3mfld}
Let $n = 3$. Let
$\JJ_i^0$ (for $i = 1,2$) be symmetric patterns of 
closed convex (or quasiconvex) sets in hyperbolic
space $\Hyp^{3} = \Hyp$ such that the limit sets
of $\JJ_i^0$ are connected.
Then any proper bijection $\phi$ between
$\JJ_1^0$ and $\JJ_2^0$ is induced by a hyperbolic isometry.
\end{theorem}

\noindent {\bf Proof:} Since limit sets are connected, we may assume
by the Scott core Theorem \cite{scott-cc}
that each $J_{1i} \in \JJ_1^0$ is the (Gromov compactified) universal
cover of a compact hyperbolic 3-manifold with incompressible
boundary. In particular,
its limit set $L_{1i}$ shares a boundary circle $C_{1i}$ with
the unbounded component of its complement. Adjoining all
the bounded components of $S^2_\infty \setminus L_{1i}$ to $L_{1i}$
we obtain $2$-disks $D_{1i}$ invariant under $\mu_l^{-1} \circ \mu_n$  
as in Theorem \ref{main} (a). Again, by Brouwer's fixed point
Theorem $\mu_l^{-1} \circ \mu_n$ has fixed points in $D_{1i}$. The
rest of the proof is as in Theorem \ref{main} (a). $\Box$

\smallskip

We next generalize Theorem \ref{main} (a) to symmetric patterns of codimension one
closed convex (or quasiconvex) sets with connected limit sets such that their stabilizers
are duality groups. This is similar to Theorem
\ref{3mfld} above. 

\begin{theorem} \label{codim1}
Let $n \geq 3$. Suppose 
$\JJ_i^0$ (for $i = 1,2$) are symmetric patterns of 
closed convex (or quasiconvex) sets in hyperbolic
space $\Hyp^{n} = \Hyp$ such that the limit sets
of $\JJ_i^0$ are connected of dimension $(n-2)$ and assume that the stabilizers
of elements of $\JJ_i^0$ (freely indecomposable codimension one
quasiconvex subgroups of $G$ by the restriction on limit sets) are duality groups.
Then any proper bijection $\phi$ between
$\JJ_1^0$ and $\JJ_2^0$ is induced by a hyperbolic isometry.
\end{theorem}

\noindent {\bf Proof:} Since limit 
sets $\LL_i^0$ of $\JJ_i^0$ have codimension one,
it follows that their stabilizers are codimension one in the big
group $G$ (the 
group acting on $\Hyp$ cocompactly). 

The  argument in this paragraph
is similar to an argument of Kapovich
and Kleiner \cite{kap-kl-pd}.
Let $G_1$ denote a stabilizer of
(some) $L \in \LL_1^0$. Since $G_1$ is a duality group, it follows
that elements of $\LL_i^0$
have the same homology as a wedge of $(n - 2)$-spheres.
By Alexander
duality, each component of the domain of discontinuity
 (= the complement of the limit set )
$\Omega (G_1) = S^{n-1}_\infty \setminus \bigcup_{L \in \LL_i^0} L$
 is acyclic.  Since $G_1$ is quasiconvex
(and hence convex-cocompact), there are only
finitely many $G_1$-orbits of such components and the 
stabilizers $H_i$, $i= 1 \cdots k$ of such components
act on them cocompactly. Therefore each $H_i$ is a $P D(n - 1)$-group.

Since 
each $H_i$ is a $P D(n - 1)$-group, the limit set of each $H_i$ is an
$(n-2)$ homology sphere $S_i$ by Theorem \ref{bestvina0}. By Alexander duality again, $S_i$
separates $S^n_\infty$ into two acyclic components (so the domain
of discontinuity of $H_i$ has two components). Adjoining
either of  these to $S_i$ gives an absolute retract (AR). 

Since  Lefschetz fixed-point Theorem holds for AR's (Remark \ref{bestvina1})
the proof of Theorem \ref{main} (a) goes through as before.
$\Box$

\subsection{Local Homology and PD(2k+1) Subgroups}
Bestvina \cite{bestvina-mmj}
 shows that Gromov boundaries of Poincar\'e
duality ($PD(m)$) hyperbolic groups are homology spheres
(Theorem \ref{bestvina0}). Thus,
if one knew some homology analogues of
properties (2), (3) in Remark \ref{genlzn} above for such spaces, Pattern Rigidity
would follow for subgroups which are $PD(2k+1)$.

We  connect the work we have done so far in this
paper to local homology properties of boundaries of hyperbolic groups
and classical techniques in algebraic topology and fixed-point theory.

\medskip

\noindent {\bf Homotopies and Coarse Topology}
 Lemma \ref{homotopy} goes through for topological manifolds and more generally,
ANR's. But more importantly for us, it generalizes to the
coarse category, where the coarse topology used
is that of Schwartz \cite{schwartz_pihes},
Farb-Schwartz \cite{farb-schwartz}, as refined and generalized
by Kapovich-Kleiner \cite{kap-kl-pd}. To see this,
first recall the following
consequence of a 
Theorem of Bestvina-Mess \cite{bestvina-mess}.

\begin{theorem}\label{bes-mess}
For a PD(n) hyperbolic group acting properly and cocompactly 
on a proper finite dimensional
simplicial complex  
$X$ with metric inherited from the simplicial structure,
there exists a compact exhaustion by compact sets $B_n$
such that the natural inclusion map of $X \setminus B_{n+1}$
into $X \setminus B_{n}$ induce isomorphisms on homology.
\end{theorem}

We shall also be using the following Theorem 
which is a result that follows from work of
Bestvina-Mess \cite{bestvina-mess} and Bestvina
\cite{bestvina-mmj} (See also Swenson \cite{swenson-mmj}, Bowditch \cite{bowditch-cutpts} and Swarup
\cite{swarup-cutpts}). 

\begin{theorem} \label{bestvina} Let $G$ be a  PD(n) hyperbolic
group
acting freely, properly, cocompactly on a contractible complex $X$
then  $H_n^{LF}(X) = H_{n-1} \partial G$. 
\end{theorem}

\begin{rmk} {\rm The isomorphism $I$ (say) of Theorem \ref{bestvina} is functorial with respect to quasi-isometries,
i.e. if $f$ is a simplicial quasi-isometry of  $X$ to itself, and $q = \partial f$ is the induced map from 
$\partial G$ to itself, then $q_\ast \circ I = I \circ f_\ast$.} \label{funct} \end{rmk}

\medskip

\noindent {\bf Approximating Quasi-isometries by Lipschitz and smooth maps}\\  Let $X$ be a convex 
contractible manifold (possibly with boundary)
of pinched negative curvature equipped with a 
 cocompact $G$-action (in particular, $G$ is Gromov-hyperbolic). Triangulating $X/G$ and lifting the triangulation to $X$, 
we have a a proper finite dimensional, locally-finite
simplicial complex  structure on
$X$ 
equipped with a proper simplicial
 cocompact $G$-action. Let $f$ be a $(K, \epsilon )$- quasi-isometry
of   $X$. Let $f^0$ be the restriction of $f$ to the zero-skeleton of the triangulation. Let $v_0, \cdots , v_k$
be vertices of a top dimensional simplex $\Delta \subset X$.
Let $d_i$ be the distance function from $v_i$.  Let $\alpha_0, \cdots , \alpha_k$ be the barycentric co-ordinates
of a point $x \in \Delta \subset X$. We define
(following Kleiner \cite{kleiner-lstr}) $\hat{f} (x)$ to be the unique point in $X$ minimizing $\sum \alpha_i d_i^2$. It follows from
work of Kleiner \cite{kleiner-lstr} that $\hat{f}$ is a Lipschitz self-map of $X$ uniformly close to $f$. Hence
$\hat{f}$ is also a  quasi-isometry with quasi-isometry constants depending on $(K, \epsilon )$ and the pinching constants
of $X$. Since $X$ is itself smooth, we can further approximate $\hat{f}$ by a smooth Lipschitz self-map of $X$ uniformly close to $f$.
For the purposes of this paper, $X$ will typically be a closed $\epsilon$-neighborhood 
$\overline{N_\epsilon (CH(\Lambda ))}$ of the convex hull of a limit
set $\Lambda$ of a convex cocompact $G$. Thus we can approximate each $\hat{f_n}$ arbitrarily closely
by a smooth map homotopic to $f_n$. Since our concern 
in this paper will be with convex hulls of limit sets of quasiconvex
subgroups in rank one symmetric spaces, we can therefore assume
that each $f_n$ is a  smooth map. {\it Note that approximating by Lipschitz maps is a considerably less delicate issue
than approximating by bi-Lipschitz maps (cf. \cite{whyte-bl}).}

Now consider a sequence $f_n$ of uniform quasi-isometries of (the vertex set
of the Cayley graph) of $G$ acting freely and cocompactly by isometries on a convex 
contractible manifold $X$ (possibly with boundary)
of pinched negative curvature. By the above discussion, we may approximate $f_n$ by 
smooth uniformly Lipschitz uniform self-quasi-isometries.

Then the coarse version of  Lemma \ref{homotopy} is

\begin{lemma} \label{coarse} Given 
\begin{enumerate}
\item a proper finite dimensional, locally-finite
simplicial complex  structure on a smooth convex manifold 
$X$ (possibly with boundary) of pinched negative curvature  with
auxiliary simplicial metric inherited from the simplicial structure, 
equipped with a proper simplicial
 cocompact action
by a Poincar\'e duality Gromov-hyperbolic group $G$, \\
\item a sequence
$f_n$ of simplicial {\bf uniform} $(K, \epsilon )$ 
quasi-isometries  of $X$ converging
uniformly on compact sets to the identity,\\
\end{enumerate}

There exist 
a positive integer $ N$ such that for all    $n \geq N$ and all $k \geq 0$,
 $f_n$ induces the identity on the locally finite homology
$H_k^{LF} (X)$.
\end{lemma}

\noindent {\bf Proof:} By the discussion preceding this Lemma,
we can assume, without loss of generality, that each $f_n$
can be approximated by smooth maps. Further, if $f_n$ converges
uniformly on compact sets to the identity, it follows (from the barycentric simplex construction of Kleiner
 \cite{kleiner-lstr} outlined above)
that so do the smooth approximants. Hence without loss of generality we may assume that $f_n$'s are smooth uniform self
quasi-isometries of $X$.

Also, since $f_n$'s are {\bf uniform} $( K, \epsilon )$ 
quasi-isometries, there exist $A_1 \geq A_2 \geq 10$ (say)
and a positive
integer $N$, such that for 
$n \geq N$ large enough, no point outside a ball of radius
$A_1$ is mapped inside a ball of radius $A_2$ under $f_n$.
Further (taking $N$ larger if necessary), we may assume by 
Lemma \ref{homotopy}, 
(since each $f_n$ is sufficiently close to the identity 
map on the ball $B_{10}$ of radius $10$ about a fixed origin $0$)
that $f_n$ restricted to $B_{10}$ is homotopic to the identity
with small tracks. By using a homotopy on a slightly smaller ball
(say of radius $9$, say) we may assume that each $f_n$ {\bf is }
the identity on $B_9$, and using straightening homotopies, 
we may also assume that no point of the complement $B_9^c$ gets mapped
to $B_9$ under $f_n$. But then the degree of the
map induced by $f_n$
on locally finite homology $H_k^{LF} (X)$ is the same as the degree
of the map $f_\ast : H_k (X, X \setminus A_1) \rightarrow  H_k (X, X \setminus A_2)$, by
Theorem \ref{bes-mess}. This is the same as the local degree of $f_n$
(see \cite{hatcher-bk} for the local degree formula)
at $0$, which in turn is $1$. Note that 
if $G$ is a hyperbolic Poincar\'e duality group of dimension $m$, $H_k (X, X \setminus A_1)$ vanishes
for $k \neq m$ and $A_1$ sufficiently large.
 This proves the Lemma. $\Box$

Now if $q_n = \partial f_n$ 
is a sequence of boundary values of {\it uniform} quasi-isometries $f_n$, 
such that $q_n$ converges uniformly to the identity map, then we may assume
that there is a point $0$ such that each $f_n$ moves $0$ through a
uniformly bounded amount. By composing with bounded track homotopies
if necessary, we may homotope $f_n$ to maps which satisfy the hypotheses 
of Lemma \ref{coarse}. 

\begin{cor} \label{alexander} Let $L = \partial X$ be the boundary of
a PD(n) hyperbolic group. 

\begin{enumerate}
\item Then $L$ is
 a compact \underline{homology manifold} with the singular (co)homology groups
of a sphere $S^d$. 
\item Let $q_i$ be a uniformly
Cauchy sequence of homeomorphisms of  $L = \partial X$
(i.e. for all $\epsilon > 0$
there exists $N$ such that for all $x \in L$, 
$d(q_m(x), q_k(x)) < \epsilon$ for all $m, k \geq N$) 
 induced by (uniform) $K, \epsilon$
 quasi-isometries $f_i$
of $X$ such that (for a fixed base-point $o$) $f_i (o)$ lies
in a uniformly bounded neighborhood
of $o$. Then there exists a positive integer
$N $ such that for all $i, j \geq N$, $q_i$ and $q_j$ induce the 
\underline{same} isomorphism
on homology groups of $L$. 
\end{enumerate}
\end{cor}

\noindent {\bf Proof:}  Assertion (1) of Corollary 
follows from Theorem \ref{bestvina0}. We provide some details here. One of the main results of
\cite{bestvina-mess} asserts that the (reduced) \v{C}ech cohomology groups of $L$
vanish except in  dimension $(n-1)$.  Bestvina \cite{bestvina-mmj} also shows that the 
(reduced) Steenrod homology groups of $L$ vanish except in  dimension $(n-1)$. Since 
$L$ is compact metrizable, Steenrod homology coincides with \v{C}ech homology (see for instance \cite{milnor-steenrod}).
Further, for locally connected metrizable compacta such as $L$, the \v{C}ech (co)-homology groups coincide
with singular (co)-homology groups (see pg. 107 of \cite{lefschetz_bk2}). Hence both
the singular (as well as \v{C}ech) homology and cohomology of $L$ coincide with that of a sphere of dimension $(n-1)$.

Using the functoriality of Remark \ref{funct}, the second assertion of the 
 Corollary now follows from the first assertion and Lemma  \ref{coarse}. $\Box$

\medskip

\noindent {\bf Scheme:} Our strategy to extend the techniques of
Theorem \ref{main} (b) beyond spheres to Poincar\'e duality
 PD(2n+1) groups (to ensure
 even dimensional boundary) is as follows: \\ 
1) Recall a consequence of an
 old  Theorem of  Lefschetz  \cite{lefschetz1},
\cite{lefschetz2}, \cite{lefschetz_bk} p.324 (for what Lefschetz calls quasicomplexes that partly generalize ANR's) generalized 
by Thompson \cite{thompson_wcx} (to weak semicomplexes that embrace
quasicomplexes, ANR's and homology manifolds in the sense of
Wilder \cite{wilder}) as
 also  Knill, \cite{knill}
Corollary 4.3, (in the general context of what Knill calls
Q-simplicial complexes) that in modern
terminology says that
 the Lefschetz fixed point Theorem holds  for 
generalized co(homological) manifolds (in the sense of
Wilder \cite{wilder}). \\
2) Use  Theorem  \ref{bestvina0} (or Corollary \ref{alexander} Assertion 1) due to
 Bestvina 
that the boundary of a hyperbolic PD(m)  group over the integers
is a homological manifold (in fact a homology
sphere) with locally connected boundary. \\
3) Finally  use Corollary \ref{alexander} Assertion 2 to conclude that
the homeomorphisms of the 
homological manifolds we have,
moving points through very small distances,
induce the identity map on homology. \\

We shall be needing the following           

\begin{theorem} (Lefschetz  \cite{lefschetz1},
\cite{lefschetz1}, \cite{lefschetz_bk} p.324, Thompson \cite{thompson_wcx}
 and Corollary 4.3 of Knill \cite{knill} ) \\
  If $Y$ is a  compact locally connected generalized homology manifold
then for any continuous map 
$f \co Y  \rightarrow Y$, 
if the Lefschetz number $A(f) \neq 0$, then $f(y) = y$ for some
$y \in Y$.
\label{knill-lfpt}
\end{theorem}

Combining Theorem \ref{knill-lfpt} with Assertion (1) of Corollary \ref{alexander}, we get 

\begin{cor} 
  If $Y$ is the boundary of a $PD(m)$ Gromov-hyperbolic group,
then for any continuous map 
$f \co Y  \rightarrow Y$ 
if the Lefschetz number $A(f) \neq 0$, then $f(y) = y$ for some
$y \in Y$.
\label{bdy-lfpt}
\end{cor}

Corollary \ref{bdy-lfpt} and  Corollary \ref{alexander} combine to give the following.

\begin{prop} \label{alexander2} Let $L = \partial X$ be the boundary of
a PD(n) hyperbolic group. 
 Let $q_i$ be a uniformly
Cauchy sequence of homeomorphisms of  $L = \partial X$
 induced by (uniform) $K, \epsilon$
 quasi-isometries $f_i$
of $X$ such that (for a fixed base-point $o$) $f_i (o)$ lies
in a uniformly bounded neighborhood
of $o$. Then there exists a positive integer
$N $ such that for all $i, j \geq N$, $q_i^{-1} \circ q_j: L \rightarrow L$ as also $q_i: L \rightarrow L$
has a fixed point. 
\end{prop}

\noindent {\bf Proof:} The second assertion of
Corollary \ref{alexander} shows that the there exists a positive integer
$N $ such that for all $i \geq N$, the 
maps $q_i$
 induce the identity map on the singular homology group.

Using the pairing between top dimensional singular
homology and cohomology, the Lefschetz number of $q_i$ or $q_i^{-1} \circ q_j$
(for all $i, j \geq N$) computed via singular cohomology 
(or equivalently, via \v{C}ech cohomology) is also the Euler characteristic. 
Corollary \ref{bdy-lfpt} now furnishes the desired conclusion. $\Box$

\medskip

The conclusion of Theorem \ref{main} (b)
for (symmetric patterns of convex hulls of limit sets of)
 PD(2k+1) subgroups 
now follows exactly along the
lines of Theorem \ref{main} (b):

\begin{theorem} \label{even}
Let $n \geq 3$. Suppose 
$\JJ_i^0$ (for $i = 1,2$) are symmetric patterns of 
closed convex (or quasiconvex) sets in hyperbolic
space $\Hyp^{n} = \Hyp$ such that the stabilizers of limit sets
of $\JJ_i^0$ are $PD(2k+1)$ 
quasiconvex subgroups of $G$.
Then any proper bijection $\phi$ between
$\JJ_1^0$ and $\JJ_2^0$ is induced by a hyperbolic isometry.
\end{theorem}

\section{Consequences and Questions}

\subsection{Rank One Symmetric Spaces}
As explained by Schwartz in Section 8 (specifically Lemma 8.1)
of \cite{schwartz_pihes}, Lemmas \ref{zoom}, \ref{eccentric}
and \ref{zariski} generalize to Complex hyperbolic space.
So do Lemmas \ref{contdiff} and \ref{scatter} (which are
generalizations of Lemmas of Schwartz \cite{schwarz-inv} ).
Thus Theorem \ref{omni} generalizes to the following.

\begin{theorem} \label{rk1}  Suppose 
$\JJ_i$ (for $i = 1,2$) are symmetric patterns of 
closed convex (or quasiconvex) sets in complex hyperbolic space $\Hyp$ of (real) dimension $n$, $n > 2$. 
For $i= 1,2$, let $G_i$ be the corresponding
uniform lattices. 
Let $H_i \subset G_i$ be 
 infinite index quasiconvex subgroups stabilizing the limit set
of some element of $\JJ_i$ and 
satisfying one of the following
conditions: 
\begin{enumerate}
\item $H_i$ is a codimension one {\it duality group}. 
\item  $H_i$ is an odd-dimensional  Poincar\'e Duality group $PD(2k+1)$
with $2k+1 \leq n-1$.
\end{enumerate} 

Then any proper bijection between
$\JJ_1$ and $\JJ_2$ is induced by a hyperbolic isometry.
\end{theorem}

\begin{rmk} {\rm For other rank one symmetric spaces (quaternionic and Cayley hyperbolic spaces), 
any quasi-isometry is a bounded distance from an isometry
by work of Pansu \cite{pansu}. Hence Theorem \ref{rk1} goes through for all rank one symmetric spaces.} \end{rmk}

\subsection{Special Disconnected Limit Sets and Intersections}
All of what we have done so far goes through with minor modifications for
disconnected limit sets, at least one of whose components
has a stabilizer $H$ of the form (1) or (2) in Theorem \ref{rk1}
above. To see this, let us retrace the argument in Theorems \ref{main}.
There we showed that for large enough $m, n$, $\mu_n^{-1} \circ \mu_m$
preserves limit sets that are spheres. The same argument shows that
for large enough $m, n$, $\mu_n^{-1} \circ \mu_m$
preserves components of limit sets of diameter bigger than (some fixed)
$\epsilon$. Since the limit set of $H$ has components whose
stabilizers are of the form (1) or (2), the arguments for 
Theorems \ref{codim1} and \ref{even} go through to prove the
existence of fixed points for $\mu_n^{-1} \circ \mu_m$. This is enough
to show the following.

\begin{cor} \label{disconnected} Let $n \geq 3$. Suppose 
$\JJ_i$ (for $i = 1,2$) are symmetric patterns of 
closed convex (or quasiconvex) sets in hyperbolic
space $\Hyp^{n} = \Hyp$ or more generally a rank one symmetric space $\Hyp$ of (real) dimension $n$. For $i= 1,2$, let $G_i$ be the corresponding
uniform lattices 
in  $\Hyp$. 
Let $H_i \subset G_i$ be an 
 infinite index quasiconvex subgroup stabilizing the (possibly
disconnected) limit set
of some element of $\JJ_i$ and 
satisfying the condition that the limit set of $H_i$ has components whose
stabilizers $H_i^{\prime}$ (obviously containing $H_i$)
are of  one of the following
forms: 
\begin{enumerate}
\item $H_i^{\prime}$ is a codimension one {\it duality group}. 
\item  $H_i^{\prime}$ is an odd-dimensional  Poincar\'e Duality group $PD(2k+1)$
with $2k+1 \leq n-1$. 
\end{enumerate}

Then any proper bijection between
$\JJ_1$ and $\JJ_2$ is induced by a hyperbolic isometry.
\end{cor}

We next state a generalization of Theorem \ref{omni} when the intersection
of some finitely many conjugates of  $H_i \subset G_i$ is of the
form (1) or (2) above. 

\begin{cor} \label{intersection} Let $n \geq 3$. Suppose 
$\JJ_i$ (for $i = 1,2$) are symmetric patterns of 
closed convex (or quasiconvex) sets in hyperbolic
space $\Hyp^{n} = \Hyp$ or more generally a rank one symmetric space $\Hyp$ of (real) dimension $n$. For $i= 1,2$, let $G_i$ be the corresponding
uniform lattices 
in  $\Hyp$. 
Let $H_i \subset G_i$ be an 
 infinite index quasiconvex subgroup and $g_1, \cdots g_m \in G$
be finitely many elements such that  $H_i^{\prime} = \bigcap_{j = 1
\cdots m} g_j H_i g_j^{-1}$ is of  one of the following
forms: 
\begin{enumerate}
\item  $H_i^{\prime}$ is a codimension one {\it duality group}. 
\item  $H_i^{\prime}$ is an odd-dimensional  Poincar\'e Duality group $PD(2k+1)$
with $2k+1 \leq n-1$.
\end{enumerate} 

Then any proper bijection between
$\JJ_1$ and $\JJ_2$ is induced by a hyperbolic isometry.
\end{cor}

\noindent {\bf Sketch of Proof:} The condition 
$H_i^{\prime} = \bigcap_{j = 1
\cdots m} g_j H_i g_j^{-1}$  implies (by  Theorems
of Short \cite{short} and Gitik-Mitra-Rips-Sageev \cite{GMRS} )
that $H_i^{\prime}$ is quasiconvex and 
that 
$\Lambda_i^{\prime} = \bigcap_{j = 1\cdots m} g_j \Lambda_i$, 
where $\Lambda_i^{\prime}$ (resp. $\Lambda_i$) represents
the limit sets of $H_i^{\prime}$ (resp. $H_i$). Since the maps
$\mu_n^{-1} \circ \mu_m$ preserve limit sets and hence their intersections
it follows that the collection of translates of joins of limit sets
of $H_i^{\prime}$ is a symmetric pattern preserved by 
$\mu_n^{-1} \circ \mu_m$. The rest of the argument proving pattern 
rigidity
is as in Theorem \ref{main}. $\Box$

\subsection{Quasi-isometric Rigidity}

Let $\GG$ be  a graph of groups  with Bass-Serre tree of spaces
$X \rightarrow T$. Let $G = \pi_1 \GG$.

(We refer the reader to \cite{msw2} specifically for the following notions:
\begin{enumerate}
\item  Depth zero raft. 
\item  Crossing graph condition. 
\item  Coarse finite type and coarse dimension.
\item  Finite depth.)
\end{enumerate}

Combining Theorems
1.5, 1.6 of \cite{msw2}  with the Pattern Rigidity theorem
\ref{omni} we have the following QI-rigidity Theorem along the lines
of Theorem 7.1 of \cite{msw2}.\\

\begin{theorem} \label{qirig}
Let $\GG$ be a finite, irreducible graph of groups such that for
the associated
Bass-Serre tree $T$ of spaces  
 no depth zero raft of $T$ is a line.
Further suppose that the vertex groups are fundamental
groups of compact hyperbolic $n$-manifolds, or more generally
uniform lattices in rank one symmetric spaces
of dimension $n$, $n \geq 3$, and 
 edge groups are all of exactly one the following types:\\
a) A  {\it duality group} of codimension one
in the adjacent vertex groups.  In this case we require in addition 
that the crossing graph condition  of Theorems 1.5, 1.6
of  \cite{msw2} be satisfied
and that $\GG$ is of finite depth.\\
b) An odd-dimensional  Poincar\'e Duality group $PD(2k+1)$ with $2k+1 \leq
 n-1$.\\
 If $H$ is a finitely generated group quasi-isometric to $G = \pi_1 ( \GG )$ then $H$ splits as
a graph $\GG^{\prime}$ of groups whose depth zero vertex groups are commensurable to those of $\GG$ and
whose edge groups and positive depth vertex groups are respectively quasi-isometric to
groups of type (a),  (b).
\end{theorem}

\noindent {\bf Proof:} By the restrictions on the vertex 
and edge groups, it automatically follows that
 all vertex and edge groups are PD groups of coarse finite type.
In Case (b), $\GG$ is
 automatically finite depth, because an infinite index subgroup of
a $PD(n)$ groups has coarse dimension at most $n - 1$. Also
the crossing graph is empty in this case hence the crossing graph 
condition  of Theorems 1.5 and 1.6 of \cite{msw2}
is automatically satisfied.

Then by  Theorems 1.5 and 1.6 of \cite{msw2}, 
$H$ splits as a graph of groups $\GG^{\prime}$
 with depth zero vertex spaces
quasi-isometric to $\Hyp = \Hyp^n$ and edge groups quasi-isometric to 
the edge groups of $\GG$ and hence respectively type (a),  (b). 
Further, the quasi-isometry respects the vertex and edge spaces of this 
splitting, and thus the quasi-actions of the vertex groups on the vertex spaces of $\GG$ preserve
the patterns of edge spaces.

By Corollary \ref{commen} the depth zero vertex groups in  $\GG^{\prime}$
 are commensurable to the corresponding groups in  $\GG$. $\Box$

Using Theorem \ref{rk1} or Corollary \ref{disconnected}, we could
get the corresponding generalizations to quasiconvex subgroups
covered by these Theorems.

\subsection{Questions}
Note that our proof of Lemma \ref{alexander} 
 does not answer the following.

\begin{question} 
Let $G$ be a PD(m) hyperbolic
group. Let $\partial G$ be its (Gromov) boundary
equipped  with a visual metric $d$. 
Does there exist $\epsilon > 0$ such that if $f$ is a homeomorphism
of $\partial G$ satisfying $d(x, f(x)) < \epsilon$ for all
$x \in \partial G$, then $f$ induces the identity map on homology?
\end{question}

Kapovich \cite{kapovich-gt} constructs convex projective
representations of the Gromov-Thurston examples; it is conceivable that these may be realized as convex cocompact
subgroups of uniform hyperbolic lattices. But there is a dearth of examples of higher
dimensional Kleinian groups in general (See \cite{kapovich-survey} for a survey.) In particular, there is a dearth
of examples in higher dimensional Kleinian groups to which Theorem \ref{rk1} applies.

There are non-ANR examples of hyperbolic Coxeter
group boundaries
coming from work of Davis \cite{davis}. These boundaries are not
locally simply connected. Doubling some of these examples
(in dimension $\geq 5$)
along their
boundaries gives the standard topological sphere $S^n$. Thus exotic (non-ENR)
homology spheres might conceivably arise as limit sets. The following
seems interesting in its own right.

\begin{question} Does there exist a convex cocompact (i.e. geometrically finite) PD(n)
hyperbolic group $G$ with non-ENR Gromov boundary
acting on $\Hyp = \Hyp^{n+1}$? Can such a $G$ appear as a codimension one quasiconvex
subgroup of a uniform lattice in $\Hyp$?
\end{question}

Fischer \cite{hfischer} has further investigated these examples. 

\begin{obs} \label{3obs} {\rm
Note that Theorem \ref{omni} combined with Corollary
\ref{disconnected} implies that pattern rigidity holds for all 
quasiconvex subgroups of hyperbolic 3-manifolds that are not virtually
free. Hence a test-case not covered by the work in this paper
 is that of symmetric patterns in
hyperbolic 3-manifolds corresponding to free quasiconvex subgroups.
This is the subject of work in progress.

Another test case is the case of symmetric patterns of
quasiconvex surface subgroups
in
hyperbolic 4-manifolds,
or, at the level of limit sets, copies of $S^1$ in $S^3$.}
\end{obs}

\begin{rmk} {\rm Much of what has been done in the context of Poincar\'e Duality groups might as well have been done in group-free language in the context of
Coarse Poincar\'e Duality spaces. (See Kapovich-Kleiner \cite{kap-kl-pd}.)}
\end{rmk}

\bibliography{patterns}
\bibliographystyle{alpha}

\bigskip

\noindent School of Mathematical Sciences, RKM Vivekananda University, \\
Belur Math, WB 711202, India 

 \noindent
\texttt{kingshook@rkmvu.ac.in, mahan@rkmvu.ac.in}

\end{document}